\documentclass{conm-p-l}

\theoremstyle{definition}

\theoremstyle{remark}

\newcommand{\abs}[1]{\lvert#1\rvert}

\begin{document}

\title[Local Densities near a Vertex]
 {Local Spectral Density and  Vacuum Energy \\
Near a Quantum Graph Vertex}

\author[S. A. Fulling]{Stephen A. Fulling}
\address{Departments of Mathematics and Physics,  Texas A\&M University,
 College Station, TX,  77843-3368}
\email{fulling@math.tamu.edu}


 \date{\today}  



\begin{abstract}
The delta interaction at a vertex generalizes the Robin
(generalized Neumann) boundary condition on an interval.  Study
of a single vertex with $N$ infinite leads suffices to determine
the localized effects of such a vertex on densities of states,
etc.  For all the standard initial-value problems, such as that
for the wave equation, the pertinent integral kernel (Green 
function) on the graph can be easily constructed from the 
corresponding elementary Green function on the real line.  From the 
results one obtains the spectral-projection kernel, local spectral 
density, and local energy density.  The energy density, which 
refers to an interpretation of the graph as the domain of a 
quantized scalar field, is a coefficient in the asymptotic 
expansion of the Green function for an elliptic problem involving 
the graph Hamiltonian; that expansion contains spectral/geometrical 
information beyond that in the much-studied heat-kernel expansion. 
\end{abstract} 

\maketitle

\section*{Introduction}

 A topic of perennial and renewed interest in quantum field theory 
is the energy of the ``vacuum'' --- that is, of the ground state of 
a field subjected to some nontrivial external condition
 \cite{Casimir,BmC,Boyer,DC,Lamo,Milton,BMM}.
 (The prototype is the electromagnetic field between two parallel 
flat conductors.)
 Although the only quantities indisputably open to experiment are
 the derivatives of the total energy with respect to parameters
 defining the configuration,
 the energy itself and even its localization in space are of 
theoretical interest, not least because energy \emph{density}
 (along with associated quantities such as pressure) acts as the 
source of the  gravitational field in general relativity 
\cite{Ford,DC}.

 From a mathematical point of view, vacuum energy
is a probe of the spectral properties of a self-adjoint 
differential operator, say $H$;
 it contains ``nonlocal'' information not extractible from the 
much-studied small-time expansion of the heat kernel 
\cite{Kirsten,Gilkey}.
It reflects the oscillatory fine structure of the eigenvalue 
distribution and, therefore, is directly related to the spectrum of 
periodic orbits \cite{Gutz,BB3,Colin,Chazarain,DG,CPS,CRR,Uribe} 
 of the nonrelativistic 
classical-mechanical 
 (or  ray-optical) system associated with our operator $H$ as 
quantum Hamiltonian (or wave operator).
 The interplay among spectral theory, dynamics, and vacuum energy 
is  fascinating and rapidly developing, with each subject gaining 
benefits from the others \cite{BmC,BD,SS,MSSvS,norman,JS}.
 Inasmuch as quantum graphs provide instructive models of spectral 
theory and semiclassical dynamics,
 they should also be communicating with vacuum energy.

Here we show that some elementary techniques recently applied 
\cite{BF} to certain traditional boundary-value problems of the 
Robin type  actually apply also to quantum graphs.
 Indeed, they may be more valuable there, because of quantum graph 
theory's supply of nontrivial problems susceptible to essentially 
one-dimensional methods.

 \section*{Infinite star graphs}

 Here we will consider only the simplest type of quantum graph, one 
consisting of a single vertex with $N$ infinite edges attached.
 (Since many of the issues we will study are basically local, there 
is a sense in which these graphs are the building blocks for all 
others.)
 The Hilbert space of the model thus consists of vector-valued
 functions,
 $u= \{u_j(x)\} \in L^2(0,\infty)^N$,
 where $x$ as the argument of $u_j$ is the distance of the point in 
question from the vertex along edge~$j$.

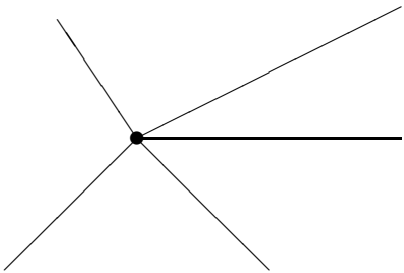
 \[ 
 \begin{picture}(150,100)(-50,-50)
 \put(0,0){\circle*{5}}
 \put(0,0){\line(1,0){100}}
 \put(0,0){\line(2,1){100}}
 \put(0,0){\line(1,-1){50}}
 \put(0,0){\line(-2,3){30}}
 \put(0,0){\line(-1,-1){50}}
 \end{picture}
   \]

 The self-adjoint operator is 
 $H = -\,\frac{d^2}{dx^2}$ with certain boundary conditions.
 We impose the usual continuity conditions,
 \begin{equation}   \label{continuity}
u_j(0) = u(0), \quad\forall j= 1,\ldots, N  .
\end{equation}
The remaining condition will be one of these:
 \begin{itemize}
 \item the Dirichlet condition, $u(0)=0\,$;
 \item the Kirchhoff, or generalized Neumann, condition,
 \begin{equation} \label{kirchhoff}
\sum_{j=1}^N u_j'(0)= 0\,;
 \end{equation}
 \item our main concern, the Exner--\v{S}eba 
or  generalized Robin condition
 \begin{equation} \label{exnerseba}
 \sum_{j=1}^N u_j'(0)= \alpha u(0), \quad\alpha>0.
 \end{equation}
 \end{itemize}
 Condition (\ref{exnerseba}) was apparently introduced in 
\cite{ES}. 
 It is often \cite{Exner,Kuchment} called the $\delta$ condition 
because it can be regarded as the effect of attaching a Dirac delta 
potential at the vertex.
At a vertex with only one edge it reduces to the so-called Robin 
(or convective cooling) condition, $u'(0)=\alpha u(0)$.
 In passing we remark that the label ``Robin'' has almost no 
historical justification \cite{GA}, but it is preferable to 
``mixed'' because ``mixed boundary conditions'' has acquired other 
meanings \cite{Kirsten,Gilkey}.

\smallskip
\textsc{Remark.} The case $\alpha <0$ can also be handled, but the 
construction of the Bondurant transform (\ref{bondurant})
 is then so different as to require a separate discussion.
 With our sign convention, $\alpha\ge0$ is the more ``physical'' 
case, where heat flows from the hotter to the cooler and $H$ has no 
negative eigenvalues.

 \section*{The Bondurant transform}  

In \cite{BF} we showed how to obtain solutions of the simplest
 problems with Robin boundary conditions from solutions of the 
corresponding Dirichlet problems.
 The term ``Dirichlet-to-Robin transform'' has sometimes been 
 misunderstood as referring to an analogue of the 
 Dirichlet-to-Neumann map, hence the justification for immediately 
naming the construction after my junior collaborator.
 Here one is \emph{not} studying the relation between 
nonhomogeneous Dirichlet data and nonhomogeneous Neumann data for a 
fixed solution;
 instead, one is constructing a \emph{new} solution to a different 
problem, with homogeneous Robin (or Exner--\v{S}eba) data replacing
Dirichlet data.

The generalization of the key formula of \cite{BF} to an infinite 
star graph with boundary conditions (\ref{continuity}) and 
(\ref{exnerseba}) is
 \begin{align}      \label{bondurant}
 u_j(x) &= (T^{-1}v)_j(x) \\
 & \equiv
  \frac1{\alpha}\left[ \frac1N \sum_{k=1}^N v_k(x) -v_j(x)\right]
 -\frac1{N^2} \int_x^\infty e^{-\alpha(s-x)/N} \sum_{k=1}^N 
v_k(s)\,ds.         \nonumber
 \end{align}

 \medskip
 \textsc{Theorem.} {\em
 If $v(t,x)$ solves a Dirichlet problem for a constant-coefficient
  partial differential equation [cf.\ (\ref{waveeq})--(\ref{cyleq})],
 then $u(t,x)$ solves the corresponding
 Exner-\v{S}eba problem (though with different initial data).
 }\medskip

\textsc {Sketch of verification.}
 (\ref{bondurant}) is obtained by solving the ordinary 
 differential equation $v=Tu$, where 
 \begin{equation}\label{bdop}
  (Tu)_j (x) \equiv \sum_{k=1}^N u_k'(x) - \alpha u_j(x),
 \end{equation}
 with the condition of decay as $x\to\infty$.
 The heuristics of finding the solution are less instructive than
 the verification that it satisfies all the required conditions.
 Since $v(0)=0$, one observes that
 \begin{itemize}
\item $u_j(0)$ is independent of $j$ (condition 
(\ref{continuity}));
\item  $Tu(x)=v(x)$, so 
 $\displaystyle\sum_{j=1}^N u_j'(0)=\alpha 
u(0)$ (condition (\ref{exnerseba}));
 \item $T$ and $T^{-1}$ commute with the partial differential 
operator, so  $u$ is still a solution.
 \hfill$\square$\hfilneg
      \end{itemize}

 \section*{Integral kernels}

 Now consider any of the following initial-value problems.
 \begin{alignat}{2}
\label{waveeq}&\mbox{Wave:} 
&\quad &u_{tt} =u_{xx}\,, \quad u(0,x) = f(x), \ u_t(0,x) = 0 \\
\label{heateq}&\mbox{Heat:}
&\quad &u_t =u_{xx}\,, \quad u(0,x) = f(x) \\
\label{quanteq}&\mbox{Quantum:}
&\quad &iu_t =-u_{xx} =Hu, \quad u(0,x) = f(x) \\
\label{cyleq}&\mbox{Cylinder:}
&\quad &u_{tt} = Hu, \quad u(0,x) = f(x), \  u(+\infty,0)= 0
 \end{alignat}
 (Of course, there are many other problems involving the operator 
$H$ on the graph that could be considered, but these form a natural 
and highly useful quartet.)
  We seek the integral kernel (Green function) that solves such a 
  problem via
 \begin{equation}\label{kernel}
 u_j(t,x) = \sum_{l=1}^N\int_0^\infty 
 G_{Sj}{}\!^l (t,x,y) f_l(y)\,dy.
 \end{equation}

 \medskip
 \textsc{Corollary.} {\em Let $G_S(t,x,y) $ be the (matrix) Green 
function for one of the initial-value problems 
 (\ref{waveeq})--(\ref{cyleq}) on the graph. 
 Let $G(t,\abs{x-y})$ be the corresponding (scalar) Green function on 
the real line (also known as the ``free'' kernel).
 Then
\begin{multline}\label{green}
 G_{Sj}{}\!^l (t,x,y) = \delta_j{}\!^l G(t,\abs{x-y}) \\
 + \left(\frac2N -\delta_j{}\!^l\right) G(t,x+y) 
 -\frac{2\alpha}{N^2}\int_x^\infty e^{-\alpha(s-x)/N} 
G(t,s+y)\,ds.
 \end{multline}
}\medskip

\textsc{Sketch of derivation.}
 The Dirichlet Green function consists of  an incident term minus 
an image term,
 \begin{align*}
 G_{Dj}{}\!^l(t,x,y) &= \delta_j{}\!^l [G(t,\abs{x-y}) - 
G(t,x+y)] \\
&\equiv G_- -G_+\,.
 \end{align*}
In operator language, we want $G_S = T^{-1} G_D T$.
 (The final $T$ is needed to get the correct initial value,
  $\delta_j{}\!^l \delta(x-y)$.) 
In kernel language, therefore, we need
\[G_S(t,x,y) = T_x^{-1} T_y^\dagger G_D(t,x,y),\]
 where $\dagger$ indicates the transpose (real adjoint).
 But 
 \begin{align}\label{transposeD}
 (T^\dagger u)_j &= -\sum_k u_k' -\alpha u_j \\
 \label{transposeT}  &= -(Tu)_j -2\alpha u_j \,.
 \end{align}
    From (\ref{transposeD}) and
$ -\partial_y G(t,\abs{x-y}) = + \partial_x G(t,\abs{x-y})$
 one sees that the incident term passes through the similarity 
transformation unchanged: $T^{-1}G_-T=G_-\,$.
        From (\ref{transposeT}) and
$-\partial_y G(t,x+y) = - \partial_x G(t,x+y)$
one has $T_y G_+ = T_x G_+$ and hence
 \begin{equation}\label{greenstep}
 G_{Sj}{}\!^l(t,x,y) =\delta_j{}\!^l [G(t,\abs{x-y}) +G(t,x+y)]
 +2\alpha T_{x,j}^{-1}[\delta_j{}\!^l G(t,x+y)] .
 \end{equation}
 (Note that the first two terms of (\ref{greenstep}) solve the true 
Neumann (not Kirchhoff) problem, $u_j'(0)=0\ \forall j$.)
 Working out the last term of (\ref{greenstep}) according to 
 (\ref{bondurant}), one obtains (\ref{green}).
\hfill$\square$\hfilneg\medskip

 The factor $\left(\frac2N -\delta_j{}\!^l\right)$ in 
(\ref{green})
 will come as no 
surprise to those who are familiar with the study of quantum graphs 
by other methods (see \cite{KS}).

\smallskip 
\textsc{Remark.} The kernel formula (\ref{green}) is correct for 
 $\alpha=0$ (the Kirchhoff condition), 
 although the intermediate steps are meaningless
 (in particular, $T^{-1}$ doesn't exist in that case).

 \smallskip
 \textsc{Main results.}
 As corollaries of the corollary, we derive formulas
 (\ref{waveker})--(\ref{wavesol}),  (\ref{specprojker}),  (\ref{specdens}),  
 (\ref{intspec}),  (\ref{cylkernel}),  (\ref{graphendens}) 
 for particular kernels and associated quantities.

\section*{The wave kernel}  

 In a one-dimensional system the simplest member of the quartet is 
the wave problem (\ref{waveeq}), for which the free Green function 
is (d'Alembert's solution)
 \begin{equation}\label{dalembert}
 G(t,z) = {\textstyle\frac12} [\delta(z-t) +\delta(z+t)].
 \end{equation}
 Applying (\ref{green}) and omitting terms that vanish for $t>0$ 
one gets
\begin{multline}\label{waveker}
 G_{Sj}{}\!^l(t,x,y) = 
 {\textstyle\frac12}\delta_j{}\!^l [\delta(x-y-t) +\delta(x-y+t)] 
\\
 {} +\frac12 \left(\frac2N -\delta_j{}\!^l\right)\delta(x+y-t) -
  \frac{\alpha}{N^2}\, e^{-\alpha(t-y-x)/N} \theta(t-y-x) ,
 \end{multline}
 where $\theta$ is the unit step function.

The meaning of (\ref{waveker}) becomes clearer when one applies
 (\ref{kernel}) to get
 \begin{multline}\label{wavesol}
  u_j(t,x)  =  {\textstyle\frac12}[f_j(x-t)+f_j(x+t)] 
-{\textstyle\frac12} f_j(t-x) \\
 {}+ \frac1N \sum_{l=1}^N f_l(t-x) 
       -\frac{\alpha}{N^2} \theta(t-x)   \int_0^{t-x}\!
  e^{-\alpha\epsilon/N}
\sum_{l=1}^N f_l(t-x-\epsilon)\,d\epsilon.
 \end{multline}
 Here we see clearly the incident wave, the immediately reflected 
and transmitted waves from a Kirchhoff vertex, 
 and some $\alpha$-dependent delayed transmission.
 The Robin case ($N=1$) is \cite{BF}
\begin{equation}\label{robsol}
 u(t,x) = 
 {\textstyle\frac12}[f(x-t)+f(x+t) + f(t-x)]
  -\alpha \theta(t-x) \int_0^{t-x} e^{-\alpha\epsilon}
  f(t-x-\epsilon)\,d\epsilon.
\end{equation}
What is the physical meaning of this time delay?
 In the context of the wave equation, the Robin or Exner--\v{S}eba 
boundary models an ideal spring, or  elastic support, to which the 
vibrating medium is attached. The spring absorbs energy from the 
medium and leaks it back out.
It is a jolly exercise in integration by parts to show that for the 
solution (\ref{robsol}) the total energy
 \begin{equation}\label{robenergy0}
E= \frac12\int_0^\infty 
 \left[\left(\frac{\partial u}{\partial t}\right)^2 +
\left(\frac{\partial u}{\partial x}\right)^2\right] dx
 +\frac{\alpha}2 u(t,0)^2
 \end{equation}
 is indeed conserved;
 the field and boundary terms individually are not (unless 
$\alpha=0$).

 For later use we note that one integration by parts in 
(\ref{robenergy0}) and use of $u'(t,0)=\alpha u(t,0)$ lead to an 
alternative formula for the total energy,
\begin{equation}\label{robenergy1/4}
 E= \frac12\int_0^\infty 
 \left[\left(\frac{\partial u}{\partial t}\right)^2 -
u\left(\frac{\partial^2 u}{\partial x^2}\right)\right] dx,
\end{equation} 
  in which the boundary term has been formally absorbed into the 
field term.

 \section*{The spectral projection kernel}

 Let $P(\lambda)$ be the orthogonal projection operator onto the 
 interval $[0,\lambda]$ of the spectral resolution of~$H$.
 Because $H$ on an infinite star graph has purely absolutely 
continuous (and nonnegative) spectrum,
  the integral kernel of $P(\lambda)$  may be written
 \[
 P(\lambda,x,y) = \int_{0}^{\sqrt{\lambda}}
 \sigma(\omega,x,y)\,d\omega,
 \]
where $\sigma$ is a well-defined matrix-valued 
function (not just a distribution in $\omega$).
 Alternatively, $\sigma$ can be defined as the
 inverse Fourier cosine transform of the wave kernel.
(It is also the 
inverse Laplace transform of the heat kernel, and so on for all 
the standard ``spectral functions'' \cite{Kirsten}.)

 That is, we write
\begin{equation}\label{fct}
G_{Sj}{}\!^l (t,x,y) =\int_0^\infty \cos (\omega t) 
 \sigma_{Sj}{}\!^l (\omega,x,y)\, d\omega
 \end{equation}
 and calculate
 \begin{align}\label{specprojker}
\sigma_{Sj}{}\!^l (\omega,x,y) &= \frac2{\pi} \int_0^\infty \cos (\omega t) 
G_{Sj}{}\!^l (t,x,y)\,dt \\ \begin{split}
 &=\frac2{\pi} \, \delta_j{}\!^l \sin(\omega x) \sin(\omega y) \\
&\quad{} +\frac{2/\pi}{\alpha^2+N^2\omega^2} \{N\omega^2 \cos[\omega(x+y)]
 +\alpha\omega \sin[\omega(x+y)]\}. 
 \end{split} \nonumber\end{align}
In a sense, (\ref{specprojker}) is the ultimate formula concerning 
the operator~$H$, since all operator functions of~$H$ can be 
calculated from it in principle and it contains all facts about 
the spectral resolution in a rather explicit form.

 Kottos and Smilansky  \cite[Sec.~3B]{KS}
 found the spectral resolution by treating the infinite star graph 
as a scattering problem.
 For each $\omega$ they provide the basis of incoming scattering 
eigenfunctions
 \begin{align}\label{scatbasis}
  \psi_j{}\!^l(x) &\equiv \delta_j{}\!^l e^{-i\omega x}
 + \left[-\delta_j{}\!^l +
 \frac1N \left(1+e^{-2i\tan^{-1}\frac{\alpha}{N\omega}} 
\right)\right] e^{i\omega x} \\
&= -2i \delta_j{}\!^l \sin(\omega x) + 2\omega\, 
\frac {N\omega -i\alpha}{ \alpha^2 +N^2\omega^2}\, e^{i\omega x} .
\nonumber \end{align}
These basis elements are orthonormal (up to the conventional factor 
$\sqrt{2\pi}$)
 and therefore
\begin{equation}\label{scatprojker} 
 \frac1{2\pi} \sum_{l=1}^N \psi_j{}\!^l(x)
  \psi_{j'}{}\!^l(y)^*
 = \sigma_{Sj}{}\!^{j'}(\omega,x,y).
 \end{equation}
A calculation verifies that (\ref{specprojker}) and 
(\ref{scatprojker}) agree. 

 \smallskip
 \textsc{Remark:} Direct construction of eigenfunctions by applying 
the Bondurant transform to orthonormal eigenfunctions of the 
Dirichlet problem, while possible, is not recommended. 
 In the present problem the 
immediate results are not orthogonal, much less normalized.
 If $\alpha=0$ they are not even linearly independent
 (because $T$ is not invertible), and one 
basis element needs to be found by a separate argument.

 \section*{The local spectral density}
 \smallskip

 Special interest attaches to the diagonal values of $\sigma$,
 \begin{multline} \label{specdens}
 \sigma_{Sj}{}\!^j(\omega,x,x) = \frac1{\pi} + 
 \frac1{\pi} \left(\frac2N -1\right) \cos(2\omega x) \\
{}  +\frac{2\alpha/\pi}{\alpha^2+N^2\omega^2}
  \left[\omega\sin(2\omega x) 
  -\frac{\alpha}N\, \cos(2\omega x)\right]. 
 \end{multline}
Clearly here the $\frac1{\pi}$ is the universal Weyl term for a 
one-dimensional system,
  the next term is the spectral effect of a Kirchhoff vertex,
and the last term is the Exner--\v{S}eba correction.
In the limit  $\alpha\to+\infty$ there is some cancellation between 
the second and third terms, resulting in
\begin{equation}\label{specdensdir} 
  \frac1{\pi}- \frac1{\pi} \cos(2\omega x) 
 =\sigma_{Dj}{}\!^j(\omega,x,x) .
 \end{equation}
as expected for a Dirichlet vertex.

Because the spectrum is continuous, 
one can't integrate (\ref{specdens}) to get a density of states.
 However, subtracting off the Weyl term and paying due attention to 
distributional limits, 
 one can obtain a meaningful global spectral density:
 \begin{align}\label{intspec}
\Delta\rho(\omega) &\equiv \int_0^\infty \sum_{j=1}^N 
 \left[ \sigma_{Sj}{}\!^j(\omega,x,x) -\frac1{\pi}  \right] dx \\
& = \left(\frac12 -\frac N4\right)\delta(\omega) 
 +\left[ \frac{N\alpha/\pi}{\alpha^2 +N^2\omega^2} 
 -\frac12\, \delta(\omega)\right] .
 \nonumber\end{align}
This expression approximates the incremental effect that such a 
vertex would have in a problem with discrete spectrum.
 The first term in (\ref{intspec}) is the Kirchhoff term;
 it vanishes when $N=2$, because a Kirchhoff vertex with exactly 
two edges is vacuous.
 The other term is the Exner--\v{S}eba correction;
 it vanishes when $\alpha=0$ because its first term 
distributionally approaches $\frac12 \delta(\omega)$ in that limit.
 Alternatively, (\ref{intspec}) can be simplified to
 \begin{equation}\label{intspecdir}
\Delta\rho(\omega)=
 -\,\frac N4\,\delta( \omega) + \frac{N\alpha/\pi}{\alpha^2 
+N^2\omega^2}\,;
 \end{equation}
 here the first term is the correct formula for a Dirichlet vertex 
and the remaining term is $O(\alpha^{-1})$ as $\alpha\to \infty$.

 \smallskip
 \textsc{Remark:} The meaning of a Dirac delta distribution in 
formulas such as (\ref{intspec}) and (\ref{intspecdir}) is that 
 the spectral counting function $N(\omega)$ has a nonzero limit as 
$\omega$ approaches $0$ from above
 ($N(\omega)$ being understood to be $0$ for negative~$\omega$).
 For example, (\ref{intspecdir}) is simply the derivative of
 the formula
 \begin{equation}\label{staircase}
  \Delta N(\omega) =  \left[-\,\frac N4 + \frac1{\pi} \,\tan^{-1} 
\frac{N\omega}{\alpha} \right]\theta(\omega)
 \end{equation}
 for the incremental effect of the vertex on the total number of 
eigenvalues in the interval $0\le\lambda\le \omega^2$.
 \smallskip

 The Bondurant method cannot be applied directly to a finite 
interval, because no transformation $T$ will work for both boundary 
conditions simultaneously. However, the appropriate operators 
 $T^{-1}$ for the two boundaries can be applied alternately to 
construct a solution as an infinite series (a sum over 
closed classical paths, generalizing the classic method of images).
  In \cite{BF} the wave kernel and hence the local and global 
spectral densities were obtained in that way for a finite interval 
with one Robin and one Dirichlet endpoint. 
 Numerical evaluation reveals the correct eigenvalues for the 
problem emerging as spikes in the global density (the counterpart 
of (\ref{intspec})). 
See \cite{SPSUS} for a related study in two dimensions.
 It should be straightforward to extend this analysis (and also the 
study of vacuum energy, etc.)\ to an arbitrary \emph{finite} star 
graph, and in principle to more complicated quantum graphs.
 It is noteworthy that in these systems no semiclassical (or 
stationary-phase) approximation is needed to obtain the 
representation of the spectrum in terms of classical paths;
 the only approximation involved is the truncation of the sum at 
some maximum path length if and when one resorts to numerics.

\section*{Heat and quantum kernels}

 The same machine (\ref{green}) can be used to treat the problems
 (\ref{heateq}) and (\ref{quanteq}), for which the free kernels
 are
 \[
 G(t,z) =  (4\pi  t)^{-1/2} e^{-z^2/4t}, \quad
  G(t,z) =  (4\pi i t)^{-1/2} e^{-z^2/4it}, 
 \]
respectively.
 (The results are qualitatively similar to (\ref{cylkernel}) below, 
with the complementary error function appearing instead of the 
exponential integral function.)
 Studying  the heat kernel is the traditional route to 
 equations like (\ref{staircase})
 for partial differential operators.
  In one-dimensional systems such as quantum graphs, however,
 the wave kernel built from (\ref{dalembert}) appears to be easier to 
 calculate.

 Using the heat kernel, the Robin case of (\ref{intspecdir})
 was obtained in \cite[Secs.\ 3.3 and 5.5]{BF}.
 That analysis extends to flat boundaries with constant $\alpha$ in 
any dimension.
The results  fit naturally with those of \cite{BB} in dimension $3$
 and \cite{SPSUS} in dimension~$2$, where the leading orders in 
boundary curvature are also included.
 All these formulas are exact in $\alpha$; of course, when the 
 heat-kernel formulas are expanded in power series in $\alpha$ they 
match and extend the relevant terms tabulated in such references 
as \cite{Kirsten,Gilkey}.
 (See also \cite{Dowker,BFSV}.)
Apart from a unifying point of view, it is not claimed that these 
results are particularly new;
 in fact, the Robin heat kernels in dimensions $1$ and $3$ were 
found in 1891 by a closely related method \cite{Bryan,Bryan3}.

 \section*{Vacuum energy density}

First, a paragraph which pure mathematicians are free to ignore:
   From a \emph{physical} point of view, 
 vacuum energy involves a relativistic (usually massless and
bosonic) field.  
 (It is of no relevance, therefore, to quantum graphs if they are 
regarded solely as models of nonrelativistic electrons in networks 
of wires.)
 Formally, the total energy corresponding to the wave operator~$H$ 
is
 \begin{equation}\label{totenergy}
 E = \frac12 \sum_{n=1}^\infty \omega_n
 =\frac12\int_0^\infty \omega\, \rho(\omega)\, d\omega
 \end{equation}
 for an operator with purely point spectrum.
 (To avoid irrelevant complications, let us also always assume that 
$H$ has no negative spectrum.)
Equally formally, the local energy density is
 \begin{equation}\label{energydens}
 T_{00}(x)=  \frac12 \int_0^\infty \omega\, \sigma(\omega,x,x)\, 
d\omega 
 \end{equation}
 (without the requirement of point spectrum).
 The origin of these expressions is the ``second-quantized'' theory
 of a field satisfying the wave equation (\ref{waveeq}), in which 
each normal mode of the field (with  frequency~$\omega$)
 becomes a quantum harmonic oscillator  (with ground state energy 
 $\frac12\omega$).
 Then (\ref{energydens}) results from the integrand of
 (\ref{robenergy1/4}), and (\ref{totenergy})  comes from 
integrating (\ref{energydens}) over all space or just from adding 
up the energies of all the modes.
Both integrals, (\ref{totenergy}) and (\ref{energydens}),
  are divergent at the upper limit and are to be 
defined by a renormalization procedure.

 Our claim is that vacuum energy should be of mathematical interest 
even in models whose physical relevance is questionable.
 Therefore, we provide here a precise
 \emph{mathematical} definition, which incorporates a particular 
renormalization prescription (whose physical rationale need not 
concern us):
{\em   Consider the cylinder kernel
 (the Green function of (\ref{cyleq}))
 on diagonal ($y=x$, $l=j$), 
find its asymptotic expansion as $t\to0$, 
 and extract the coefficient of the term 
proportional to $t$, times~$-\frac12$; 
 this is the vacuum energy density, $T_{00}(x)$.
 When appropriate, integrate over $x$ 
 (and sum over $j$ in our graph problem),
before taking $t$ to $0$, to 
 define a total  energy,~$E$.}

 The intuition behind this definition is the following.
 Let $G(t,x,y)$ be the cylinder kernel of the problem under study.
 (For an infinite star graph it is the matrix~$G_S\,$.)
 Then
 \begin{equation}\label{cyldiag}
 G(t,x,x) = \int_0^\infty e^{-\omega t} \sigma(\omega,x,x)\, 
d\omega,
\end{equation}
and (when ``appropriate'') its trace is
 \begin{equation}\label{cyltrace} 
 \mathop{\mathrm{Tr}}\nolimits  G(t)
\equiv \int_0^\infty \sum_{j=1}^N G_j{}\!^j (t,x,x)\, dx
 = \sum_{n=1}^\infty e^{-\omega_n t}.
 \end{equation}
Now take $-\frac12 \,\frac{\partial}{\partial t}$ of 
 (\ref{cyldiag}) and (\ref{cyltrace}) and let $t$ approach~$0$,
 formally obtaining (\ref{energydens}) and (\ref{totenergy}), 
respectively.
 Systematically throwing away the terms of negative order in the 
Laurent expansions (and a possible logarithmic term), one arrives 
at our definition.

To find the vacuum energy of an infinite star graph, 
 one can apply the Bondurant machine one 
more time, to the free cylinder kernel
 \begin{equation}\label{freecyl}
  G(t,z) = \frac{t/\pi}{t^2+z^2}\,,
 \end{equation}
 obtaining
 \begin{multline}\label{cylkernel}
  G_{Sj}{}\!^l = \delta_j{}\!^l \, \frac{t/\pi}{t^2+(x-y)^2}
 +\left(\frac 2N -\delta_j{}\!^l\right)\frac{t/\pi}{t^2+(x+y)^2} 
\\
 + \frac{2\alpha}{\pi N^2} \, e^{\alpha(x+y)/N}  
 \mathop{\mathrm{Im}}\nolimits
 \left[e^{-i\alpha t/N} \mathop{\mathrm{Ei}}\nolimits
  \left( \frac{i\alpha t}{N} 
 -\frac{\alpha}N\,(x+y) \right)\right].
 \end{multline}
 It follows that
\begin{equation}\label{graphendens}
  T_{00}(x) =
 \left(1-\frac 2N\right) \frac1{8\pi x^2} 
 + \frac{\alpha}{2\pi N^2 x} +\frac{\alpha^2}{\pi N^3}\,
 e^{2\alpha x/N} 
 \mathop{\mathrm{Ei}}\nolimits
 \left(-\,\frac{2\alpha x}{N}\right) . 
 \end{equation}

 The most important parameter in this problem is the dimensionless 
product $\alpha x$.
 Therefore,
 at short distance the Kirchhoff term dominates:
 \begin{equation}\label{nearlim}
 T_{00}(x) \sim  \left(1-\frac 2N\right) \frac1{8\pi x^2}
  + \frac{\alpha}{2\pi N^2 x} 
 +\frac{\alpha^2}{\pi N^3} \,\ln|\alpha x| +O(\alpha^2x^2),
 \end{equation}
whereas at large distance the energy density is almost pure Dirichlet:
 \begin{equation}\label{farlim}
T_{00}(x) \sim \frac 1{8\pi x^2} +O(\alpha^{-1} x^{-3}).
 \end{equation}
The nonintegrable $O(x^{-2})$ singularity in (\ref{nearlim})
would interfere with calculating a total energy by integration 
over~$x$,  even if the edges were finite.  The renormalization 
procedure implicit in our definition does not commute with the 
spatial integration, however, and it leads to a finite total 
energy \cite{BGH,systemat,norman}.
 As agreed, we will not delve here into the physical issues thereby 
raised (which are still somewhat controversial).

 In a sense the calculation based on the cylinder kernel was 
unnecessary, given the spectral formulas (\ref{specprojker}),
 (\ref{specdens}), and (\ref{intspec}).
 Indeed, (\ref{graphendens}) can be obtained from (\ref{cyldiag})
 and (\ref{specdens}) and a Laurent expansion,
 or even immediately from (\ref{energydens}) using the subtracted 
spectral density appearing in the integrand of (\ref{intspec}).
 (The subtraction of the Weyl term corresponds to the subtraction 
of the leading Laurent term.  In other problems, such as higher 
dimensions, additional subtractions would be necessary.)
Similar methods were used in \cite{RS} to find the effect of a flat 
Robin boundary in any dimension.
In general, however, it is easier to find the small-$t$ expansion 
of a cylinder kernel than to find the detailed spectral resolution;
 as for heat kernels, one expects useful calculations to run in the 
opposite direction.  
 Romeo and Saharian \cite{RS} also gave complicated integral 
formulas for the vacuum energy and energy density of a finite 
interval with  two Robin boundaries.
 Our methods will instead give these quantities as infinite sums 
over classical paths \cite{FL}.

\section*{A broader perspective}

Now let $H$ be a generic differential operator
 (self-adjoint, elliptic, positive, second-order, with scalar 
principal symbol)
 in dimension~$d$.
 Let $T$ be the cylinder kernel and $K$ be the heat kernel.
 Then
\[
\begin{aligned}
\mathop{\mathrm{Tr}}\nolimits T &=  \int_0^\infty e^{-t\omega} \,dN
 =\int_0^\infty e^{-t\omega} \rho(\omega)\,d\omega, \\
\mathop{\mathrm{Tr}}\nolimits  K &=  \int_0^\infty e^{-t\lambda  } \,dN
 =\int_0^\infty e^{-t\omega^2  } \rho(\omega)\,d\omega,
 \end{aligned}
\quad
 \begin{aligned} T(t,x,x) &=  \int_0^\infty e^{-t\omega} 
 \sigma(\omega,x,x)\,d\omega, \\
 K(t,x,x) &=  \int_0^\infty e^{-t\omega^2  } 
\sigma(\omega,x,x)\,d\omega,
  \end{aligned}
 \]
where $N$ now is the number of eigenvalues less than or equal to 
$\lambda=\omega^2$. 
  
It is known that as $t\to0$ the global quantities have expansions 
of the forms
 \begin{equation}\label{cylseries}
\mathop{\mathrm{Tr}}\nolimits  T \sim
\sum_{s=0}^\infty e_s t^{-d+s}
+\sum^\infty_{\genfrac{}{}{0pt}2{\scriptstyle s=d+1}{
\scriptstyle s-d \mbox{ \scriptsize odd}}} f_s t^{-d+s} \ln 
t,
 \end{equation}
  \begin{equation}\label{heatseries}
\mathop{\mathrm{Tr}}\nolimits
  K \sim \sum_{s=0}^\infty b_s t^{(-d+s)/2}.
 \end{equation}
The local quantities, $T(t,x,x)$ and $K(t,x,x)$, have (nonuniform 
in~$x$) expansions of precisely the same respective forms;
 we do not introduce a separate notation for their coefficients.
{\em
 If $d-s$ is even or positive,
 \begin{equation}\label{evenrel}
e_s= \pi^{-1/2} 2^{d-s} \Gamma((d-s+1)/2) b_s\,.
\end{equation}
 If $d-s$ is odd and negative, then
 \begin{equation}\label{oddrel}
f_{s} = \frac{(-1)^{(s-d+1)/2}2^{d-s+1}  }{
\sqrt{\pi}\, \Gamma((s-d+1)/2)} \, b_s\,,
 \end{equation} 
 but, most strikingly, in that case
 \begin{equation}\label{undetermined}
\mbox{$e_{s}$ is undetermined by the $b_{r}$}\,.
 \end{equation}
}

I have expounded in detail elsewhere  \cite{FG,systemat,norman}
 how the connections (\ref{evenrel})--(\ref{undetermined}) come 
about:
  The $b_s$ are proportional to coefficients in the high-frequency 
asymptotics of the Riesz means \cite{Hardy,Hor}
 of  $N$ (or of $\int\!\sigma$)
  with respect to $\lambda$, while
 the $e_s$ and $f_s$ are proportional to coefficients in the
asymptotics of the Riesz means  with respect to $\omega$.
 Alternatively, they are related through the poles of the zeta 
functions of $H$ and $\sqrt{H}$ \cite{Gilkey4,GG,BM}.

 The new half of the cylinder-kernel coefficients
 (those in (\ref{undetermined}))
---  of which the first, $e_{d+1}\,$, is proportional to the vacuum 
  energy ---
 are a new set of moments of the spectral distribution.
 {\em What are they good for, mathematically?\/}
 Unlike the old ones, they are \emph{nonlocal} in their dependence 
on the geometry of the domain and the coefficient functions of~$H$.
  They probe, in a comparatively crude, averaged way,
  the oscillatory spectral structures that are correlated more 
  precisely with 
 periodic and closed classical paths.
  Thus they embody, at least partially, the global dynamical structure 
of the system;
 they are a half-way house between the heat-kernel coefficients
 and a full semiclassical closed-orbit analysis.

\end{document}